\documentclass{elsarticle}
\usepackage{amssymb,theorem,amsmath}
\newtheorem{thm}{Theorem}[section]
\newtheorem{lem}{Lemma}[section]
\newproof{pf}{\bf Proof}

{\theorembodyfont{\normalfont}
\newtheorem{rmk}{Remark}[section]}
\newtheorem{cor}{Corollary}[section]
\newdefinition{defn}{Definition}[section]
{\theorembodyfont{\normalfont}
\newtheorem{exmp}{Example}[section]}

\makeatletter

\let \al=\alpha
\let \be=\beta
\let \var=\varphi
\let \vare=\varepsilon

\let \de=\delta
\let \th=\theta
\let \la=\lambda

\let \q=\quad

\let \med=\medskip
\let \smal=\smallskip
\let \dps =\displaystyle

\newcommand{\R}{\mathbb{R}}

\def\C{{\rm |\kern
-4.6pt{\rm C}}}
\def\N{{\rm I\kern
-4.0pt{\rm N}}}

\begin{document}
\begin{frontmatter}

\title{Global Dynamics for  Lotka-Volterra Systems with  Infinite Delay and Patch Structure}
\author{Teresa Faria\footnote{ 
 Fax:+351 21 795 4288; Tel: +351
21 790 4929.}}

\address{Departamento de Matem\'atica and CMAF, Faculdade de Ci\^encias, Universidade de Lisboa\\ Campo Grande, 1749-016 Lisboa, Portugal\\
tfaria@ptmat.fc.ul.pt}



\begin{abstract} We study some aspects of the global dynamics of an $n$-dimensional Lotka-Volterra system with infinite delay and patch structure, such as extinction, persistence,  existence and global attractivity of a positive equilibrium. Both the cases of an irreducible and reducible linear community matrix are considered, and no restriction on the signs of the intra- and inter-specific delayed terms is imposed. Although the system is not cooperative, our approach often uses comparison results applied to an auxiliary cooperative system.
Some models in recent literature are  generalised, and results improved.
\end{abstract}

\begin{keyword}  
Lotka-Volterra system; infinite delay; patch-structure;  global attractivity; persistence; extinction.\\
{\it 2010 Mathematics Subject Classification}:  34K20, 34K25, 34K12, 92D25
\end{keyword}

\end{frontmatter}

\section{Introduction} 
\setcounter{equation}{0}
In recent years, mathematicians and biologists have been analysing biological models  given by differential equations  with  time-delays and  patch-structure. Models with patch-structure are frequently quite realistic, to account for heterogeneous environments and other biological features,  where   single or multiple species  are distributed over several different patches or classes, with migration among them. Time-delays are very often present  in models from population dynamics, neurosciences, ecology, epidemiology, chemistry and other sciences. Moreover,  infinite delays have been considered in equations used in  population dynamics since the works of Volterra, to translate the cumulative effect of the past history of a system. Typically,  the ``memory functions" appear as integral kernels and,  although defined in the entire past, the delay  should be introduced in such a way that its effect diminishes when going back in time.  

\med

In this paper, 
 the following patch-structured Lotka-Volterra system with both infinite distributed and discrete delays is considered: 
\begin{equation}\label{eq0}
\begin{array}{ll}
\displaystyle{x_i'(t)=x_i(t)\biggl(b_i-\mu_i x_i(t)-\sum_{j=1}^n a_{ij} \int_0^{\infty}K_{ij}(s)x_j(t-s)ds\biggr)} \\
\displaystyle{\hskip 1.5cm +\sum_{j\ne i,j=1}^n (\vare_{ij}\al_{ij}x_j(t-\tau_{ij})-\al_{ji}x_i(t)), \quad i=1,2,\ldots,n}.
\end{array}
\end{equation}
Here, $\mu_i> 0$, $b_i,a_{ij} \in\R,$  and, for $i\ne j$, $ \al_{ij} \geq 0,\tau_{ij}\ge 0$, $\vare_{ij}\in (0,1]$, $ i,j=1,\dots,n$; the kernels $K_{ij}: [0,\infty) \to [0,\infty)$ are $L^1$ functions, normalized so that 
\begin{equation}\label{1.2}
\int_0^{\infty} K_{ij}(s)\, ds=1, \quad \mbox{for} \ i,j=1, \dots,n.
\end{equation}
Moreover, we suppose that for all $i$ the linear operators defined by $L_{ii}(\var)=\int_0^{\infty} K_{ii}(s)\var (-s)\, ds$, for $\var:(-\infty,0]\to\R$ bounded, are non-atomic at zero, which amounts to have $K_{ii}(0)=K_{ii}(0^+)$. 

 System \eqref{eq0} serves as a population model  for the growth of single or multiple species distributed over $n$ patches or classes: $x_i(t)$ is the density of the population on patch $i$, with $b_i$ and $\mu_i$ as its  usual Malthusian growth rate and (instantaneous) self-limitation coefficient, respectively, $a_{ii}$ and $a_{ij}\, (i\ne j)$ are respectively the intra- and inter-specific delayed acting coefficients;    $\al_{ij}\, (i\ne j)$ are the dispersal rates of populations moving from patch $j$ to patch $i$, and $\tau_{ij}$  the times taken during this dispersion; the coefficients $\vare_{ij}\in (0,1]$ appear to account for some lost of the populations during migration from one patch to another. Frequently, one takes 
$\varepsilon_{ij}=e^{-\gamma_{ij} \tau_{ij}}$ for some $\gamma_{ij}> 0,\ i,j=1,\dots,n,i\ne j$, cf. e.g. \cite{Takeuchi:2006}.
Denoting
$$d_{ij}:=\vare_{ij}\al_{ij}\ {\rm for}\ i\ne j,\q 
\be_i:=b_i-\sum_{j\ne i} \al_{ji},
$$
\eqref{eq0} is written as
\begin{equation}\label{eq}
\begin{array}{ll}
\displaystyle{x_i'(t)=x_i(t)\biggl(\be_i-\mu_i x_i(t)-\sum_{j=1}^n a_{ij} \int_0^{\infty}K_{ij}(s)x_j(t-s)ds\biggr)} \\
\displaystyle{\hskip 2cm +\sum_{j\ne i,j=1}^n d_{ij}x_j(t-\tau_{ij}), \quad i=1,2,\ldots,n},
\end{array}
\end{equation}
where $\be_i\in\R, \mu_i>0,a_{ij}\in\R, d_{ij}\ge 0, \tau_{ij}\ge 0$ and the kernels $K_{ij}$ are as above.

With model \eqref{eq0} in mind, in the present paper  some aspects of  the asymptotic behaviour of solutions to delayed Lotka-Volterra systems \eqref{eq} are analysed.
Although not very meaningful in biological  terms,  actually all the techniques and results in this paper apply to  systems with several bounded delays or even infinite delays in the migration terms, which leads to more general systems of the form
\begin{equation}\label{eq1}
\begin{array}{ll}
\displaystyle{x_i'(t)=x_i(t)\biggl(\be_i-\mu_i x_i(t)-\sum_{j=1}^n a_{ij} \int_0^{\infty}K_{ij}(s)x_j(t-s)ds\biggr)} \\
\displaystyle{\hskip 1cm +\sum_{j\ne i,j=1}^n \sum_{p=1}^m d_{ij}^{(p)}x_j(t-\tau_{ij}^{(p)}), \q i=1,2,\ldots,n},
\end{array}
\end{equation}
with $d_{ij}^{(p)},\tau_{ij}^{(p)}\ge 0$, or
\begin{equation}\label{eq2}
\begin{array}{ll}
\displaystyle{x_i'(t)=x_i(t)\biggl(\be_i-\mu_i x_i(t)-\sum_{j=1}^n a_{ij} \int_0^{\infty}K_{ij}(s)x_j(t-s)ds\biggr)} \\
\displaystyle{\hskip 1cm +\sum_{j\ne i,j=1}^n d_{ij}\int_0^{\infty}G_{ij}(s)x_j(t-s)ds, \quad i=1,2,\ldots,n},
\end{array}
\end{equation}
with the kernels $G_{ij}\ge 0$   being $L^1$ functions with $L^1$-norm one. Moreover, our method can be easily addapted to  Lotka-Volterra systems with   continuous coefficients and discrete delays depending on $t$. 


Due to the biological interpretation of the model,  only positive or non-negative solutions should be considered admissible. On the other hand, there are natural constraints on admissible phase spaces  for functional differential equations (FDEs) with infinite delay (cf.~Section 2). To deal with such kind of equations, a careful choice of a so-called `fading memory space' as phase  space is in order, see e.g. \cite{HaleKato,HMN},  and  one must consider  {\it bounded} initial conditions. Thus, our framework accounts 
only for solutions of \eqref{eq} with  initial conditions of the form
\begin{equation}\label{initial}
x_i(\theta)=\varphi_i(\theta), \ \theta \in (-\infty,0], \quad \varphi_i(0)>0, \ i=1,\ldots,n,
\end{equation}
where $\varphi_i$ are non-negative and bounded continuous functions on $(-\infty,0]$. 
\med
  
There is an immense   literature on FDEs of Lotka-Volterra type, 
and it is impossible to mention all the relevant contributions. 
  The present investigation was motivated by several papers, among them those of Takeuchi et al.~\cite {Tak06a,Takeuchi:2006}, Liu \cite{Liu}, and Faria  \cite{Faria:2010,Faria:2013a}.  For other related papers, we  refer to    \cite{DingHan,Faria:2013b,Muroya:2004,Wang}, also for further references.
  
   In \cite{Liu}, Liu considered a cooperative model for a species  following a delayed logistic law, with the  population structured in several classes and no delays in the migration terms, of the form
\begin{equation}\label{Liu}
\begin{array}{ll}
\displaystyle{x_i'(t)=x_i(t)\Big [b_i-\mu_{i} x_i(t)+\sum_{p=1}^mc_i^{(p)}x_i(t-\sigma_i^{(p)})\Big ]+\sum_{j=1}^n d_{ij}x_j(t), \ i=1,\ldots,n},
\end{array}
\end{equation}
where  $\mu_i>0,b_i>0$ and  $c_i^{(p)},d_{ij},\sigma_i^{(p)}\ge  0$ for $i,j=1,\dots,n, p=1,\dots,m$. Moreover, in \cite{Liu}  only the case  $D=[d_{ij}]$ an  irreducible matrix was studied, and the further quite restrictive conditions  $(b_i+\sum_{j=1}^n d_{ij})/(\mu_i-\sum_{p=1}^mc_i^{(p)})=k$ for $\, 1\le i\le n$ ($k$ a positive constant) were imposed. On the other hand,  Takeuchi et al. \cite{Takeuchi:2006}  studied the  system
\begin{equation}\label{Tak}
\begin{array}{ll}
\displaystyle{x_i'(t)=x_i(t)\big (b_i-\mu_i x_i(t)\big )+\sum_{j\ne i, j=1}^n \big(e^{-\gamma_{ij} \tau_{ij}}\al_{ij}x_j(t-\tau_{ij})-\alpha_{ji}x_i(t)\big),\ i=1,\ldots,n}, 
\end{array}
\end{equation}
where $\mu_i>0,b_i\in\R$ and $ \al_{ij},\tau_{ij},\gamma_{ij}\ge  0$ for $i,j=1,\dots,n$, $j\ne i$. Note that \eqref{eq0} is a natural generalization of \eqref{Tak}, obtained by the addition of  interacting terms with infinite delay.
Again,  only the case of an irreducible matrix $D=[d_{ij}]$, where now $d_{ij}=e^{-\gamma_{ij} \tau_{ij}}\al_{ij}$,  was studied in \cite{Takeuchi:2006}. 

In  \cite{Faria:2013a}, the author analyses several aspects of the asymptotic behaviour of solutions to the more  general cooperative system 
\begin{equation}\label{Far}
\begin{array}{ll}
\dps x_i'(t)=x_i(t)\Big [b_i-\mu_ix_i(t)+ \sum_{p=1}^m c_i^{(p)}x_i(t-\sigma_i^{(p)})\Big ]+\sum _{ j=1}^n \sum_{p=1}^m d_{ij}^{(p)}x_j(t-\tau_{ij}^{(p)}),\ i=1,\dots, n,\
\end{array}
\end{equation}
where: $ b_i\in\R, \mu_i>0$ and $ c_i^{(p)}, d_{ij}^{(p)}, \sigma_i^{(p)}, \tau_{ij}^{(p)}\ge 0 $, for $i,j=1,\dots,n,\, p=1,\dots,m$. The situations of  $D=[d_{ij}]$ an irreducible or a reducible  matrix were both addressed. Note that models \eqref{Liu} and \eqref{Tak} are particular cases of \eqref{Far}.
\med

For the patch structured Lotka-Volterra models \eqref{eq}, \eqref{eq1} or  \eqref{eq2}, for simplicity we write
 $A=[a_{ij}]$ and $ D=[d_{ij}]$, with $d_{ii}:=0$, and  where $d_{ij}:=\sum_{p=1}^m d_{ij}^{(p)}$ for  \eqref{eq1}. The matrix $M(0)=diag\, (\be_1,\dots,\be_n)-D$, i.e.,
\begin{equation}\label{M(0)}
\displaystyle{M(0)=\left[
\begin{array}{cccc}
\be_1& d_{12}&\cdots&d_{1n} \\
d_{21}& \be_2&\cdots&d_{2n} \\
\vdots&\vdots&\ddots&\vdots \\
d_{n1}&d_{n2}&\cdots&\be_n \\
\end{array}
\right]},
\end{equation}
may be interpreted as the \textit{linear community matrix}. For coefficients $a_{ij}\in\R$, we shall use the standard notation 
$$a_{ij}^-=\max (0,-a_{ij}),\q a_{ij}^+=\max (0,a_{ij}).$$
 Throughout the paper, together with $M(0)$ we shall consider the matrix
\begin{equation}\label{1.9}
N_0=diag\, (\mu_1,\dots,\mu_n)-[a_{ij}^-].
\end{equation}
The algebraic properties of $M(0)$ and $N_0$ will play a crucial role in the global dynamics of the system.
\med

For simplicity, this paper deals with system \eqref{eq}, rather than \eqref{eq1} or \eqref{eq2}, and addresses its  global asymptotic behaviour, in what concerns its dissipativity and persistence,  extinction of the populations, and the existence and global attractivity of a positive equilibrium.
As in  the cited papers \cite{Liu,Takeuchi:2006}, most papers dealing with patch structured models only analyse the situation of an irreducible linear community matrix. Here,  both the cases of $M(0)$ irreducible and reducible are  considered. Of course, if $M(0)$ is irreducible, sharper criteria can be obtained, namely a threshold criterion of exchanging of global attractivity between the trivial solution and a positive equilibrium. Rather than Lyapunov functional techniques, the approach exploited here is based on  comparison results and monotone techniques (see \cite{Smith}) applied to an auxiliary cooperative system, coupled  with  theory of M-matrices.

\med

The contents  of the paper are now  briefly described. Section 2 is a preliminary section, where an abstract formulation to deal with \eqref{eq} is set, and some notation  and auxiliary results  are given,
 including some  known properties from matrix theory; also we prove some important  estimates used throughout the paper. Section 3 provides criteria for  the local stability and global attractivity of the trivial equilibrium -- in biological terms, the latter  translates as the extinction of the populations in all patches. In Section 4, we consider  the particular case of  \eqref{eq} with all coefficients $a_{ij}\le 0$, thus a cooperative Lotka-Volterra system, and investigate its persistence  and  global asymptotic convergence to an equilibrium. Finally, Section 5 is  devoted to the study of the persistence, 
the existence and the global attractivity  of a positive equilibrium for the general model \eqref{eq}.

 \section{Preliminaries: abstract framework, notation and auxiliary results}
\setcounter{equation}{0}

In this preliminary section, we first set an abstract framework to deal with  system \eqref{eq}.  In view of the  unbounded delays, the problem  must be  carefully formulated by defining an appropriate Banach phase space where the problem is well-posed. 

 Let $g$ be a
function satisfying the following properties:
\med

 {(g1)} $g:(-\infty ,0]\to [1,\infty)$ is a non-increasing continuous function, $g(0)=1$;

 {(g2)} $\displaystyle{\lim_{u\to 0^-}{{g(s+u)}\over {g(s)}}=1}$ uniformly on $(-\infty ,0]$;

 {(g3)} $g(s)\to \infty$ as $s\to -\infty$.

\med
\noindent 
For $n\in\N$, define  the Banach space $UC_g=UC_g(\R^n):=\big\{ \phi\in C((-\infty, 0];\R^{n}) : \sup_{s\le 0}{{|\phi(s)|}\over {g(s)}}<\infty,
{{\phi(s)}\over {g(s)}}\ {\rm is\ uniformly\ continuous\ on}\ (-\infty, 0]\big\},$
with the norm
$$\|\phi\|_g=\sup_{s\le 0}{{|\phi(s)|}\over {g(s)}},$$
where $|\cdot|$ is a chosen norm  in $\R^n$. 
Consider also  the space $BC=BC(\R^n)$ of bounded  continuous functions $\phi:(-\infty, 0]\to
\R^n$.
It is clear that $BC\subset UC_g$.
%

The space $UC_g$ is an admissible phase space for $n$-dimensional FDEs with infinite delay
(cf. \cite{HaleKato, HMN}) written in the abstract form
\begin{equation}\label{2.1}
\dot x(t)=f(t,x_t),
\end{equation}
where $f:D\subset \R\times UC_g\to \R^n$ is continuous and, 
as usual, segments of solutions in the phase space $UC_g$ are denoted by $x_t$, $x_t(s)=x(t+s), s\le 0$. When $f$ is regular enough and the initial conditions are \textit{bounded}, it is known that 
the initial value problem is well-posed, in the sense that  there exists a unique solution $x(t)$ of  the problem $\dot x(t)=f(t,x_t), x_\sigma=\var\in BC$, denoted by $x(t;\sigma,\var)$ in $\R^n$ or  $x_t(\sigma,\var)$ in $UC_g$; for autonomous systems $\dot x(t)=f(x_t)$ and $\var\in BC$, the solution of $\dot x(t)=f(x_t),\, x_0=\var$ is simply denoted by $x(t;\var)\in\R^n$ and $x_t(\var)\in UC_g$. Moreover,  bounded positive orbits of \eqref{2.1} are precompact in $UC_g$ \cite{HaleKato,HMN}.

\med

An appropriate formulation for problem \eqref{eq}-\eqref{initial} is set as follows.
From Lemma 4.1 in  \cite{Faria:2011}, for any $\de>0$ there is a continuous function $g$ satisfying (g1)--(g3) and such that
\begin{equation}\label{2.2}
\int_0^{\infty} g(-s)K_{ij}(s)\, ds<1+\de,\q i,j=1,\dots,n.
\end{equation}
Whenever it is necessary, one fixes  a positive $\de$ and inserts the  problem into the phase space $UC_g$, where $g$ is any function satisfying the above conditions (g1)-(g3) and \eqref{2.2}. Of course, if one considers the more general system \eqref{eq2}, one should demand that $g$ also satisfies the conditions $\int_0^{\infty} g(-s)K_{ij}(s)\, ds<1+\de,\ i,j=1,\dots,n.$

In the space $UC_g$, a vector $c$  is identified with the constant function $\psi(s)=c$ for $s\le 0$. A vector $c$  in $\R^n$ is said to be {\it positive} (respectively {\it non-negative}) if all its components are positive (respectively non-negative).
We use the notation  $\R^n_+=\{x\in\R^n:x\ge 0\}$, and $BC^+=BC^+(\R^{n})=\{ (\var,\psi)\in BC: \var(s),\psi(s)\ge 0$ for all $s\le 0\}$.  In view of the biological meaning of  \eqref{eq}, 
the framework is restricted  to   positive or  non-negative initial conditions.
As a  set of admissible initial conditions for \eqref{eq}, we take the  subset $BC_0^+$ of $BC^+$, $BC^+_0=\{ (\var,\psi)\in BC^+:\var(0)>0,\psi(0)>0\}$. It is easy to see that all the coordinates of solutions with initial conditions in $BC^+$, respectively  $BC_0^+$, remain non-negative, respectively positive, for all $t\ge 0$ whenever they are defined (see e.g. \cite{Smith}).

A system \eqref{2.1} is said be {\it cooperative} if it satisfies the quasi-monotonocity condition in p.~78 of Smith's monograph \cite{Smith}: whenever $\var,\psi \in BC^+,\var\le \psi$ and $\var_i(0)=\psi_i(0)$ holds, then $f_i(\var)\le f_i(\psi)$, for $1\le i\le n$. Hence,
  system \eqref{eq} is cooperative if and only if $a_{ij}\le 0$ for all $i,j=1,\dots,n$.

The following crucial estimates will be used throughout the paper. For  similar arguments, cf.~\cite{Faria:2013b}.
 
\begin{lem}\label{lem1} 
Let $x:\R\to\R^n$ be a continuous function with $x_0=\var\in BC$. For each $j\in\{ 1,\dots,n\}$, suppose that there are constants $M\in\R$ and $t_0>0$ such that $x_j(t)\le M$, respectively $x_j(t)\ge M$, for $t\ge t_0$. Then, for any $\vare>0$  there exists $T_0\ge t_0$ such that
$$
\displaystyle{ \int_0^{\infty}K_{ij}(s)x_j(t-s)\, ds\le M+\vare,\q t\ge T_0,\,  1\le i\le n,}$$
respectively
$$
\displaystyle{ \int_0^{\infty}K_{ij}(s)x_j(t-s)\, ds\ge M-\vare,\q t\ge T_0,\, 1\le i\le n.}$$
\end{lem}

\begin{pf} Suppose that  $x_j(t)\le M$ for $t\ge t_0$.
 Fix $\vare >0$, and take $g$ satisfying conditions (g1)-(g3) and (2.2) with $0<\de (1+\de) \le \vare $. Since $x_j(t)$ is bounded from above on $(-\infty,\infty)$, take  $K>0$ such that $\sup_{t\in\R}x_j(t)  \le K$, and choose $T> 0$ such that $K/g(-T)<\de$.
For $t\ge T_0:=T+t_0$, from \eqref{1.2} and \eqref{2.2} we have
$$
 \begin{array}{ll}
\dps{\int_0^{\infty}K_{ij}(s)x_j(t-s)\, ds}&=\dps{\int_0^{T}K_{ij}(s)x_j(t-s)ds+\int_T^{\infty}K_{ij}(s)x_j(t-s)ds}\\
&\dps{\le M\int_0^{T}K_{ij}(s)\, ds+\int_T^{\infty}g(-s)K_{ij}(s)\frac{K}{g(-s)}ds}\\
 &\dps{\le M\int_0^{T}K_{ij}(s)\, ds+\int_T^{\infty}g(-s)K_{ij}(s)\frac{K}{g(-T)}ds}\\
 &\dps{\le M\int_0^{T}K_{ij}(s)\, ds+\de\int_T^{\infty}g(-s)K_{ij}(s)\, ds}\\
 &\dps{\le M+\de (1+\de)\le M+\vare.}
 \end{array}
 $$
The other inequality is proven in a similar way.\end{pf}

We now recall some notation and results from matrix theory. An $n\times n$ matrix $M=[m_{ij}]$ is said to be {\it cooperative} if all its  off-diagonal entries are nonnegative: $m_{ij}\ge 0$ for $i\ne j$. 
 
 Let  $\sigma (M)$ be the spectrum of $M$.  The {\it spectral bound}  $s(M)$ of $M$ is defined as
  $$s(M)=\max\{ Re\, \la: \la  \in\sigma (M)\}.$$ 
It is well-know that if $M$ is cooperative and irreducible, then $s(M)\in\sigma (M)$ and there is a positive  eigenvector associated with $s(M)$ (see e.g.~\cite{BP,Fiedler}). 

A square matrix  $M$ is said to be an {\it M-matrix}  (respectively  {\it non-singular M-matrix}) if  all its  off-diagonal  entries  are non-positive and all its
eigenvalues  have a non-negative (respectively positive) real part. 

\begin{lem}\label{lem2} For a square matrix $M=[m_{ij}]$ with $m_{ij}\le 0$ for $i\ne j$, the following conditions are equivalent:\vskip 0cm
(i) $M$ is a non-singular M-matrix;\vskip 0cm
(ii) $M$ is  an M-matrix and is non-singular;\vskip 0cm
(iii) all principal minors of $M$ are positive;\vskip 0cm
(iv) there is a positive vector $v$ such that $Mv>0$;\vskip 0cm
(v) $M$ is  non-singular and $M^{-1}\ge 0$.

\end{lem}

There are many other equivalent ways of defining   non-singular M-matrices, as well as M-matrices, see  \cite{BP,Fiedler} for a proof of Lemma \ref{lem2} and further properties of these matrices.  In \cite{Fiedler}, non-singular M-matrices and  M-matrices are also designated by {\it matrices of classes} $K$ and $K_0$, respectively. The notation is far from being  uniform, and  many authors call {\it M-matrices}  the matrices defined here as {\it non-singular M-matrices}.  

From these definitions, it is apparent that for a cooperative matrix $M$,  $s(M)\le 0$ (respectively $s(M)<0$) if and only if
$-M$ is an M-matrix (respectively a non-singular M-matrix).

For \eqref{eq},  the matrix $M(0)$  defined in \eqref{M(0)} is cooperative.
If $D=[d_{ij}]$ is irreducible, then $M(0)$ is  irreducible as well.   For cooperative   and irreducible matrices,  the following lemma is useful.

\begin{lem}\label{lem3} If $M=[m_{ij}]$ is a cooperative   and irreducible matrix, then $s(M)>0$ if and only if there exists a positive vector $v$ such that $Mv>0$.
\end{lem}

\begin{pf} Let $M=[m_{ij}]$ be cooperative and irreducible. Then $s(M)$ is an eigenvalue of $M$ with a positive associated eigenvalue $v$, thus  $s(M)>0$ implies that $Mv=s(M)v>0$ for some positive vector $v$ (cf.~\cite{BP}). 

Conversely,  let $v$ be a positive vector such that $Mv>0$, and for the sake of contradiction suppose that $s(M)\le 0$. Then $-M$ is an M-matrix, or, in other words, for any $\de>0$ the matrix $\de I-M$ is a non-singular M-matrix  \cite{Fiedler}; this implies that $(\de I-M)^{-1}\ge 0$. Choose $\de >0$ small so that $Mv>\de v$.  Then we have
$v=(\de I-M)^{-1}(\de I-M) v\le 0$, a contradiction.
\end{pf}

\section{Stability  of the trivial equilibrium }\label{sec2}
\setcounter{equation}{0}

The standard notions  of stability and attractivity used here are always defined in the context of the set $BC^+_0$ of admissible initial conditions, as recalled below. 

\begin{defn} An equilibrium $x^*\ge 0$   of \eqref{eq}  is said to be {\it stable} 
if for any $\vare> 0$  there is $\delta=\de (\vare)>0$ such that $\|x_t(\var)-x^*\|_g<\vare$ for all $\var \in BC^+_0$ with $\|\var-x^*\|_g<\delta$ and $t\ge 0$;  $x^*$ is said to be
{\it globally attractive} if $x(t)\to x^*$ as $t\to\infty$, for all  solutions $x(t)$ of \eqref{eq} with initial conditions $x_0=\var \in BC^+_0$; and   $x^*$  is {\it globally asymptotically stable (GAS)}  if it is stable and globally attractive. 
\end{defn}

In this section, we address the  stability and attractivity of the trivial equilibrium. When \eqref{eq} refers to a population model, the global attractivity of 0 means the extinction of the populations in all patches.

\begin{thm}\label{thm21} For system \eqref{eq},
(i) if $s(M(0))<0$, then the equilibrium 0 is hyperbolic and locally asymptotically stable;
(ii) if $s(M(0))>0$,  then 0 is unstable. \end{thm}

\begin{pf}
We have already observed that $s(M(0))<0$ if and only if $-M(0)$ is a non-singular M-matrix.

(i) Assume that  $s(M(0))<0$. The linearization of \eqref{eq} at zero is given by
\begin{equation}\label{lin0}
x_i'(t)=\be_ix_i(t)+\sum_{j\ne i}d_{ij}x_j(t-\tau_{ij}), \quad i=1,2,\ldots,n.
\end{equation}
Denote $B=diag\, (\be_1,\dots, \be_n)$.
The characteristic equation for \eqref{lin0} is
\begin{equation}\label{char0}
\det \Delta (\lambda)=0,\ \ {\rm where}\ \ \Delta (\lambda)=M(\lambda)-\lambda I
\end{equation} 
and
$$
\displaystyle{M(\lambda)=\left[
\begin{array}{cccc}
\be_1& d_{12}e^{-\lambda \tau_{12}}&\cdots&d_{1n} e^{-\lambda \tau_{1n}}\\
d_{21}e^{-\lambda \tau_{21}}& \be_2&\cdots&d_{2n} e^{-\lambda \tau_{2n}}\\
\vdots&\vdots&\ddots&\vdots \\
d_{n1}e^{-\lambda \tau_{n1}}& d_{n2}e^{-\la \tau_{n2}}&\cdots&\be_n \\
\end{array}
\right]=:B+D(\la)}.
$$
Since $-M(0)$ is a non-singular M-matrix,
from a result in \cite{Faria:2008}  (which can be generalised to linear FDEs with infinite delay) $x=0$ is asymptotically stable as a solution of \eqref{lin0}, for all values of the delays.

 (ii) Assume now  that $s(M(0))>0$.
 
 First consider the case of $D=D(0)$ an irreducible matrix. 
Observe that $\Delta(0)=M(0)$ and $\Delta(\la_1)> \Delta(\la_2)$ for $0\le \la_1<\la_2$.
Since the matrices $\Delta(\la), 0\le \la < \infty$, are irreducible and cooperative, then $s(\Delta(\la))\in \sigma (\Delta(\la))$, and $\la \mapsto s(\Delta(\la))$ is continuous and  strictly decreasing on $[0,\infty)$. Clearly, $s(\Delta(\la))\to -\infty$ as $t\to\infty$; together with the condition $s(\Delta(0))>0$, this implies the existence of a unique $\la^*>0$ such that $s(\Delta(\la^*))=0$. But $s(\Delta(\la^*))\in \sigma (\Delta(\la^*))$, or, in other words, $\la^*$ is a characteristic root for \eqref{lin0}.  This proves that \eqref{lin0} is unstable, hence 0 is unstable as a solution of \eqref{eq}.

Next, consider the case  of $D$ reducible. After a simultaneous permutation of rows and columns, for each $\la\ge 0$ the matrix $D(\la)$ is written in a triangular form as
$$
\displaystyle{D(\lambda)=\left[
\begin{array}{ccc}
D_{11}(\la)& \cdots&D_{1\ell}(\la)\\
&\ddots & \\
0&\cdots&D_{\ell \ell}(\la) \\
\end{array}
\right]},
$$
where $D_{lm}(\la)$ are $n_l\times n_m$ matrices, with $D_{ll}(\la)$ irreducible blocks and $\sum_{l=1}^\ell n_l=n$. Only to prove the result for $\ell=2$ is  needed, since the general case will follow by induction. 

For $\la\ge 0$, let $
D(\lambda)=\left[
\begin{array}{cc}
D_{11}(\la)&D_{12}(\la)\\
0&D_{22}(\la) \\
\end{array}
\right],$
where $D_{11}(\la),D_{22}(\la)$ are irreducible, and
 write $M(\la),\Delta(\la)$ in the form
$$\begin{array}{ll}
\displaystyle{M(\al)=B+D(\la)=\left[
\begin{array}{cc}
M_{11}(\la)&M_{12}(\la)\\
0&M_{22}(\la) \\
\end{array}
\right],}\\
\displaystyle{\Delta(\la)=M(\la)-\la I_n=\left[
\begin{array}{cc}
M_{11}(\la)-\la I_{n_1}&M_{12}(\la)\\
0&M_{22}(\la)-\la I_{n_2} \\
\end{array}
\right].}
\end{array}$$
We have $\sigma (M(0))=\sigma (M_{11}(0))\cup \sigma (M_{22}(0))$ and 
$\sigma (\Delta (\la))=\sigma (M_{11}(\la)-\la I_{n_1})\cup \sigma (M_{22}(\la)-\la I_{n_2})$. Hence,  $s(M_{ii}(0))>0$ either for $i=1$ or $i=2$, and by the irreducible case we deduce that there exists $\la^*>0$ such that  $0\in \sigma (M_{ii}(\la^*)-\la^* I_{n_i})$. This shows that $\la^*>0$ is a solution of the characteristic equation \eqref{char0}, and thus 0 is unstable.
\end{pf}

In the case of cooperative systems, we shall prove that if $M(0)$ is irreducible, then $s(M(0))>0$ is a sharp criterion for persistence, conf.~Theorem \ref{thm31}.  For the moment,  a result for the extinction of all populations is given.

\begin{thm}\label{thm22} Consider \eqref{eq}, and assume that $N_0$ is a non-singular M-matrix. Then,  all positive solutions of \eqref{eq} are   bounded. Moreover, if
there exists a positive vector $q=(q_1,\dots,q_n)$ which satisfies the conditions
 \begin{equation}\label{2.3}
 \begin{array}{ll}
\displaystyle{
\mu_i q_i-\sum_{j=1}^n c_{ij}q_j  >0,}\\
\displaystyle{ \be_iq_i+\sum_{j \neq i} d_{ij} q_j\le 0,
\ \ \ i=1,\dots,n, } \end{array}
\end{equation}
then  all
the populations go extinct
in every patch, i.e., all positive solutions $x(t)$ of \eqref{eq} satisfy $x(t)\to 0$ as $t\to\infty$.
\end{thm}

 \begin{pf} 
 Consider the  
 following auxiliary cooperative system:
\begin{equation}\label{2.4}
\begin{array}{ll}
\displaystyle{x_i'(t)=x_i(t)\biggl(\be_i-\mu_i x_i(t)+\sum_{j=1}^n a_{ij}^- \int_0^{+\infty}K_{ij}(s)x_j(t-s)ds\biggr)}, \\
\displaystyle{\hskip 2cm +\sum_{j\ne i} d_{ij}x_j(t-\tau_{ij})=:f_i(x_t), 
 \quad i=1,2,\ldots,n},
\end{array}
\end{equation}
where as before $a_{ij}^-=\max (0,-a_{ij})$,  and observe that the solutions of \eqref{eq} satisfy
$$
\begin{array}{ll}
\displaystyle{x_i'(t)\le x_i(t)\biggl(\be_i-\mu_i x_i(t)+\sum_{j=1}^n a_{ij}^- \int_0^{\infty}K_{ij}(s)x_j(t-s)\, ds\biggr)}, \\
\displaystyle{\hskip 2cm +\sum_{j\ne i} d_{ij}x_j(t-\tau_{ij}), 
 \quad i=1,2,\ldots,n.}
\end{array}
$$

By Theorem 5.1.1 of \cite{Smith}, the solution $x(t)$ of each initial value problem \eqref{eq}-\eqref{initial} satisfies $x(t)\le X(t)$, $t\ge 0$, where $X(t)$ is the solution of  \eqref{2.4}-\eqref{initial}.  Therefore, it is enough to prove the statements of the theorem for \eqref{2.4}.

Since $N_0$ is a non-singular M-matrix, there is $q=(q_1,\dots,q_n)>0$  such that $N_0q>0$, i.e.,
$\mu_i q_i-\sum_{j=1}^n a_{ij}^-q_j  >0, \, i=1,\dots,n$ (cf.~Lemma \ref{lem2}).
For $L>0$ sufficiently large, we have
\begin{equation}\label{2.5}
\displaystyle{f_i(Lq)=Lq_i\Big [\be_i-L(\mu_iq_i-\sum_{j=1}^n a_{ij}^-q_j )\Big ]+L\sum_{j \neq i} d_{ij} q_j <0, \q  i=1,\dots,n.}
\end{equation}

Consider solutions $x(t)=x(t;\var)$  of \eqref{2.4}-\eqref{initial}. From Smith's results (cf. Corollary 5.2.2 in \cite{Smith}), and since bounded positive orbits in $UC_g$ are precompact, this implies that$$x(t;\var)\le x(t;Lq)\searrow x^*\q {\rm for}\q \var\le Lq,$$
where $x^*=(x_1^*,\dots,x_n^*)$ is necessarily an equilibrium of \eqref{2.4} (recall that $\omega$-limit sets are invariant sets). In particular, this proves that all positive solutions of \eqref{2.4} are bounded.

Next, we assume \eqref{2.3} and prove that $L_i:=\limsup_{t\to\infty}\limits \frac{x_i(t)}{q_i}=0$ for all $i$.

Let $L_i=\max_j L_j$. If $L_i>0$, by the fluctuation lemma  there is a sequence $(t_k), t_k\to\infty,$ with $x_i(t_k)\to L_iq_i, x_i'(t_k)\to 0$. For $\vare>0$ small and $k$ large, the application of Lemma \ref{lem1} yields
\begin{equation}\label{2.6}
x_i'(t_k)\le x_i(t_k)[\be_i-\mu_ix_i(t_k)+\sum_j a_{ij}^- q_j (L_i+\vare)]+(L_i+\vare) \sum_{j\ne i} d_{ij}q_j.
\end{equation}
By letting $k\to\infty, \vare\to 0^+$, the above formula and \eqref{2.3} lead to
\begin{equation}\label{2.7}
\begin{array}{ll}
\displaystyle{
0\le L_i q_i\Big [\be_i-L_i(\mu_iq_i-\sum_{j=1}^n a_{ij}^-q_j )\Big ]+L_i\sum_{j \neq i} d_{ij} q_j }\\
\q =\displaystyle{
 L_i \Big [\be_iq_i+\sum_{j \neq i} d_{ij} q_j 
 -L_iq_i(\mu_iq_i-\sum_{j=1}^n a_{ij}^-q_j )\Big ]<0.}
\end{array}
\end{equation}
This is a contradiction, and the proof is complete.
\end{pf}

\begin{rmk}\label {rmk1} Condition \eqref{2.3} reads as
 \begin{equation}\label{2.8}
 M(0)q\le 0\ \ {\rm and }\ \ N_0\, q>0,
 \end{equation}
 where $M(0)$ and $N_0$ are given by  \eqref{M(0)} and  \eqref{1.9}; this is equivalent to saying that  $N_0$ is a  non-singular M-matrix and 
 $$M(0) N_0^{-1}v\le 0$$ for some positive vector $v$.
Clearly, \eqref{2.8} also implies that   $-M(0)$ is an M-matrix \cite{Fiedler}; this condition also translates as $s(M(0))\le 0$. Note however that the converse is not true: in fact, even if $M(0)$ is an irreducible matrix with $s(M(0))\le 0$ and $N_0$ is a non-singular M-matrix, then there exist positive vectors $v$ and $q$ such that
 $ M(0)v\le 0\ \ {\rm and }\ \ N_0\, q>0,$
 but one cannot conclude that there is one positive vector $q$ satisfying simultaneously $M(0)q\le 0$ and $N_0\, q>0$, and the extinction of the populations cannot be derived.  This is illustrated below by a counter-example.
 \end{rmk}
 
 \begin{exmp} Consider the system  \eqref{eq} with $n=2$, $\be_1=\be_2=-2, d_{12}=1,d_{21}=\frac{7}{2}$ and $\mu_1=-a_{12}=1, \mu_2=\frac{13}{45}, a_{21}=-\frac{1}{10}, a_{11}=a_{22}=0$:
\begin{equation}\label{2.9}
\begin{array}{ll}
\dps{x_1'(t)=x_1(t)\biggl(-2-x_1(t)+  \int_0^{\infty}K_{12}(s)x_2(t-s)\, ds\biggr)+x_2(t-\tau_1)}\\
\dps{x_2'(t)=x_2(t)\biggl(-2-\frac{13}{45}x_2(t)+\frac{1}{10} \int_0^{\infty}K_{21}(s)x_1(t-s)\, ds\biggr)+\frac{7}{2}x_1(t-\tau_2)}.
\end{array}
\end{equation}
with  delays $\tau_1,\tau_2\ge 0$ and positive kernels $K_{12},K_{21}$ satisfying \eqref{1.2}. With the previous notation,  
$$
M(0)=\left[
\begin{array}{cc}
-2&1 \\
\frac{7}{2}&-2 \\
\end{array}
\right],\ N_0=\left[
\begin{array}{cc}
1&-1 \\
-\frac{1}{10}& \frac{13}{45} \\
\end{array}
\right].$$
Clearly $N_0$ is a non-singular M-matrix and $s(M(0))<0$. The positive vectors $v=(1,v_2)$ satisfying $M(0)v\le 0$ are the ones for which $\frac{7}{4}\le v_2\le 2$; a positive vector $q=(1,q_2)$ satisfies $N_0q>0$ if and only if $\frac{9}{26}<q_2<1$. Hence there is no vector $q>0$ satisfying both conditions \eqref{2.8}. In this example, the trivial equilibrium is not a global attractor, since $(1,\frac{3}{2})$ is a positive equilibrium of \eqref{2.9}.
\end{exmp}

The next result follows clearly from the the proof of Theorem \ref{thm21}. 

\begin{thm}\label{thm23}  If
there exists a positive vector $q=(q_1,\dots,q_n)$  satisfying $M(0)q< 0$ and $N_0\, q\ge 0$,
then the equilibrium 0 of \eqref{eq} is GAS.
\end{thm}

The case of no patch structure,  has been studied by the author in   \cite{Faria:2010} (see also Faria and Oliveira \cite{Faria:2008}), where the local stability and attractivity of a positive equilibrium was investigated, but not the  extinction, for which sufficient conditions are given below.

 \begin{cor}\label{cor21} 
 Consider \eqref{eq0} with $\al_{ij}=0$ for $1\le i,j\le n$ (no patch structure):
$$
x_i'(t)=x_i(t)\bigg (b_i-\mu_i x_i(t)-\sum_{j=1}^n a_{ij} \int_0^{\infty}K_{ij}(s)x_j(t-s)\, ds\bigg ),\ i=1,\dots,n,$$
where all the coefficients and kernels are as in \eqref{eq0}. As in \eqref{1.9}, denote $N_0=diag\, (\mu_1,\dots,\mu_n)-[a_{ij}^-]$. If either (i)  $N_0$  is a non-singular M-matrix and  $b_i\le 0,\, 1\le i\le n,$ or (ii)  $N_0 q\ge 0$ for some positive vector $q$, and $b_i< 0,\, 1\le i\le n,$ then all positive solutions satisfy $x(t)\to 0$ as $t\to\infty$.\end{cor}

In the case of  competitive systems, the next corollary generalises and improves Theorem 3.3 in \cite{Faria:2013a}.

 \begin{cor}\label{cor22} Consider \eqref{eq}, with $\mu_i>a_{ii}^-$ and $a_{ij}\ge 0$ for $j\ne i, \, i,j=1,\dots,n$. If there is a positive vector $q$ such that $M(0)q\le 0$, then the equilibrium 0  is globally attractive. In particular, this holds if
 either $s(M(0))< 0$, or $s(M(0))= 0$ and $M(0)$ is irreducible.
\end{cor}

\begin{pf} In this situation, $N_0$ reads as $N_0=diag\, (\mu_1-a_{11}^-,\dots,\mu_n-a_{nn}^-)$.
The first assertion follows immediately from Theorem \ref{thm22}. If $s(M(0))< 0$, then there is a vector $q>0$ such that $M(0)q<0$. Moreover,  if  $M(0)$ is  irreducible, since it  is also cooperative, then $s(M(0))= 0$ implies the existence of a positive vector $q$ such that $M(0)q=0$ \cite{Fiedler}.
\end{pf}

\section{Persistence and  stability for the cooperative Lotka-Volterra system}\label{sec3}
\setcounter{equation}{0}

In this section,  attention is   devoted to  the  cooperative case of  \eqref{eq}, written here as
\begin{equation}\label{3.1}
\begin{array}{ll}
\displaystyle{x_i'(t)=x_i(t)\biggl(\be_i-\mu_i x_i(t)+\sum_{j=1}^n c_{ij} \int_0^{\infty}K_{ij}(s)x_j(t-s)\,ds\biggr)} \\
\displaystyle{\hskip 2cm +\sum_{j\ne i,j=1}^n d_{ij}x_j(t-\tau_{ij}), \quad i=1,2,\ldots,n},
\end{array}
\end{equation}
where $c_{ij}=-a_{ij}\ge 0$ for $i,j=1,\dots,n$. For \eqref{3.1}, the matrix $N_0$ given in \eqref{1.9} is rewritten as $N_0=diag\, (\mu_1,\dots,\mu_n)-[c_{ij}]$. 

%

For the definitions of persistence and dissipativity used below, see e.g.~\cite{Kuang,ST}.

\begin{defn} A system $x'(t)=f(x_t)$ with $S\subset BC$ as set of admissible initial conditions is said to be {\it persistent}  if any  solution $x(t;\var)$ with initial condition $\var\in S$ is bounded away from zero, i.e.,
$$\liminf_{t\to\infty} x_i(t;\var)>0,\q 1\le i\le n,$$ 
for any any $\var\in S$; and the system is said to be {\it dissipative}  if there is a positive constant $K$ such that, given any $\var\in S$, there exists $t_0=t_0(\var)$  such that 
$$|x_i(t,\var)|\le K,\q {\rm for}\q 1\le i\le n,\, t\ge t_0.$$ 
Clearly, $S=BC^+_0$ for \eqref{eq}.
\end{defn}

\begin{thm}\label{thm31} If there is a positive vector $v$ such that $M(0)v>0$, then \eqref{3.1} is  persistent;
moreover, there is a positive equilibrium. In particular, this is the case if $s(M(0))>0$ and $M(0)$ is irreducible.
\end{thm}

\begin{pf} Write \eqref{3.1} in the form
$$x_i'(t)=f_i(x_t),\q i=1,\dots,n.$$
For $v=(v_1,\dots,v_n)>0$ such that $M(0)v>0$ and $l>0$ small, we  obtain
$$
f_i(lv)=l\bigg (\be_iv_i+\sum_{j \neq i} d_{ij} v_j \bigg )-l^2v_i \bigg(\mu v_i-\sum_{j=1}^n c_{ij}v_j \bigg)>0,\  i=1,\dots,n.
$$
Hence,  there exists a positive equilibrium $x^*$ with $x(t;lv)\nearrow x^*$; moreover, since the system  \eqref{3.1} is cooperative, $x(t;\var)\ge x(t;lv)$ if $l>0$ is sufficiently small so that $\var\ge lv$  \cite{Smith}. This shows the  persistence of \eqref{3.1}.
The last assertion of the theorem follows from  Lemma \ref{lem3}.
\end{pf}

A criterion for  the global attractivity of a positive  equibrium for \eqref{3.1} is now established.

\begin{thm}\label{thm32}
Assume  there is a vector $v>0$ such that $M(0)v>0$ and  $N_0=diag\, (\mu_1,\dots,\mu_n)-[c_{ij}]$ is a non-singular M-matrix. If $x^*$ is a  positive equilibrium (whose existence is given by Theorem \ref{thm31}) and $M(0)x^*>0$, 
then $x^*$ is the unique positive equilibrium  of system \eqref{3.1} and is globally attractive. 
\end{thm}

\begin{pf} Consider  vectors $q>0, v>0$ such that $N_0q>0,M(0)v>0$.
Using the above notation, for $L>>1$ and $0<l<<1$, we have $f_i(Lq)<0$  and $f_i(lv)>0,\, 1\le i\le n$. Thus 
 there exist positive equilibria $x^*,y^*$, with
$$x(t;lv)\nearrow x^*, \ \ x(t;Lq)\searrow y^*,\ \ \ {\rm for }\ \ 0<l<<1<<L,$$
and $$x(t;lv)\le x(t;\var)\le x(t;Lq),\ \ \ {\rm for }\ \ lv\le \var\le Lq.$$
In particular, all positive solutions of  \eqref{3.1} are bounded and bounded away from zero.
To show  the uniqueness and global attractivity of a positive equilibrium, it is sufficient to prove  the following: 
\med

{\it Claim:} if $x^*=(x_1^*,\dots,x_n^*)$ is a positive equilibrium of  \eqref{3.1}, then any positive solution $x(t)=(x_1(t),\dots,x_n(t))$ of  \eqref{3.1} satisfies
$$\liminf_{t\to\infty} x_i(t)\ge x_i^*,\ \ \ i=1,\dots,n.$$

First observe that  a positive equilibrium   $x^*$ of \eqref{3.1} satisfies
\begin{equation}\label{3.3}
x_i^*\Big [\be_i-\mu_i x_i^*+\sum_j c_{ij} x_j^*\Big ]+\sum_{j\ne i} d_{ij} x_j^*=0,\q 1\le i\le n,
\end{equation}
or, in other words,
\begin{equation}\label{3.4}
M(0)x^*=x^*\otimes N_0x^*,
\end{equation}
where we use the notation  $u\otimes v=(u_1v_1,\dots,u_nv_n)$ for $u=(u_1,\dots,u_n),v=(v_1,\dots,v_n)\in\R^n$. By assumption $M(0)x^*>0$, or equivalently  $N_0x^*>0$.

To prove the claim, effect the changes of variables $\bar x_i(t)=x_i(t)/x_i^*$ in  \eqref{3.1}, and define $\ell_i:=\liminf_{t\to\infty} \bar x_i(t) >0$.

Choose $i$ such that $\ell_i=\min_j \ell_j$. 
We now drop the bars for simplicity, and consider a sequence $t_k\to \infty$ with $x_i'(t_k)\to 0$ and $x_i(t_k)\to \ell_i$. 
For any $\vare\in (0,\ell_i)$ and $k$ sufficiently large, from Lemma \ref{lem1} we get
$$x_i'(t_k)\ge x_i(t_k)\Big [ \be_i-\Big (\mu_ix_i^*x_i(t_k)-(\ell_i-\vare)\sum c_{ij}x_j^*\Big )\Big]+
(\ell_i-\vare)\frac{1}{x_i^*}\sum_{j\ne i} d_{ij} x_j^*.$$
By taking limits $k\to\infty, \vare\to 0^+$,  we obtain
$$0\ge \ell_i(1-\ell_i)\Big (\be_i+\frac{1}{x_i^*}\sum_{j\ne i} d_{ij} x_j^*\Big )= \ell_i(1-\ell_i)(N_0x^*)_i.$$
This yields $\ell_i\ge 1$, which proves the claim.
\end{pf}

\begin{cor}\label{co31}
Assume that $N_0=diag\, (\mu_1,\dots,\mu_n)-[c_{ij}]$ is a non-singular M-matrix and $\be_i>0$ for $1\le i\le n$. Then, there is a  positive equilibrium   of system \eqref{3.1}, which is a global attractor. \end{cor}

\begin{pf} If $\be_i>0$ for $1\le i\le n$, then $M(0)v>0$ for $v=(1,\dots,1)$. Moreover, if $x^*$ is a positive equilibrium of \eqref{3.1}, then by \eqref{3.4} we have $N_0 x^*>0$. The result follows from Theorem \ref{thm32}.
\end{pf}

We now treat the generalisation of model \eqref{Far}  obtained by introducing  infinite delay.

\begin{cor}\label{cor32} Consider 
\begin{equation}\label{3.5}
\begin{array}{ll}
\displaystyle{x_i'(t)=x_i(t)\biggl(\be_i-\mu_i x_i(t)+ c_i \int_0^{\infty}\! \! K_i(s)x_i(t-s)\, ds\biggr)+\sum_{j\ne i} \sum_{p=1}^m d_{ij}^{(p)}x_j(t-\tau_{ij}^{(p)}), \ i=1,\ldots,n},
\end{array}
\end{equation}
where:  $\be_i\in \R, \mu_i>0$ and $c_i,d_{ij}^{(p)},\tau_{ij}^{(p)}\ge 0$; $K_i:[0,\infty)\to [0,\infty)$ are in $L^1$ with $L^1$-norm equal to 1, $1\le i,j\le n$. Consider $M(0)$ given by \eqref{M(0)}, where $d_{ij}=\sum_{p=1}^m d_{ij}^{(p)}$.  If $\mu_i>c_i$ for $1\le i\le n$, then:\vskip 0cm
(i) if there is a positive vector $q$ such that $M(0)q\le 0$,  the equilibrium 0  is a global attractor; (ii) if there is a positive vector $q$ such that $M(0)q>0$,  there exists a positive equilibrium $x^*$ which is  a global attractor.
\end{cor}

When $M(0)$ is irreducible,  a threshold criterion for \eqref{3.5} is as follows:

\begin{cor}\label{cor33} Consider \eqref{3.5} with  $M(0)$ irreducible and $\mu_i>c_i$ for $1\le i\le n$.
Then: (i) if $s(M(0))\le 0$,  the equilibrium 0  is a global attractor; (ii) if $s(M(0))>0$,  there exists a positive equilibrium $x^*$ which is  a global attractor.
\end{cor}

\begin{rmk} System \eqref{3.5} generalizes both \eqref{Liu} and \eqref{Far}. Not only the model is more general, but also  Corollaries \ref{cor32} and \ref{cor33} provide stronger  criteria than the ones in \cite{Faria:2013a,Liu}. In fact, for   \eqref{Liu},    Liu \cite{Liu} assumed that  $[d_{ij}]$ is irreducible, $\mu_i,b_i>0$, all the other coefficients are non-negative with $\mu_i>\sum_{p=1}^mc_i^{(p)}$,   and proved that if the constants $\al_i^*:=(b_i+\sum_{j=1}^n d_{ij})/(\mu_i-\sum_{p=1}^mc_i^{(p)})$, $1\le i\le n$, are all equal to some constant $k$,    then the equilibrium $x^*=(k,\dots,k)$ is a global attractor of all positive solutions; while in Faria \cite{Faria:2013a} the existence and global attractivity of a positive equilibrium was proven simply under the  assumptions of $\mu_i>\sum_{p=1}^mc_i^{(p)}$ and $b_i+\sum_{j=1}^n d_{ij}>0,\ 1\le i\le n$.

\end{rmk}

\section{Persistence and  stability for the general Lotka-Volterra system}
\setcounter{equation}{0}

 We now return to the general case of the Lotka-Volterra model \eqref{eq} with no prescribed signs for the   interaction coefficients $a_{ij}$, whose extinction was already studied in Section 2. Sufficient conditions for   dissipativeness,  persistence, and global attractivity of a positive equilibrium will be given. 

In what follows, $M(0)$ and $N_0$ are as in \eqref{M(0)} and \eqref{1.9}. For the  case of $M(0)$ irreducible,  first observe that  there are no non-trivial equilibria on the boundary of the non-negative cone $\R^n_+$; moreover, if the system is dissipative and 0 is unstable, this implies  the existence of a positive equilibrium. A more exact result is stated in the lemma below.

\begin{lem}\label{lem41} If  $s(M(0))>0$ and  \eqref{eq} is dissipative, then \eqref{eq} has a non-trivial equilibrium $x^*\ge 0$. If in addition $M(0)$ is an irreducible matrix,  \eqref{eq}
has a positive equilibrium $x^*$.
\end{lem}

\begin{pf}
 Consider the ODE system associated with \eqref{eq}, given by
 \begin{equation}\label{4.2}
x_i'(t)=x_i(t)\biggl(\be_i-\mu_i x_i(t)-\sum_{j=1}^n a_{ij} x_j(t)\biggr)+\sum_{j\ne i}d_{ij}x_j(t), \ i=1,\dots,n. 
\end{equation}
Clearly, \eqref{eq} and \eqref{4.2} share the same equilibria.
By assumption, \eqref{4.2} is dissipative.  Since the non-negative cone $\R^{n}_+$ is forward  invariant
for \eqref{4.2}, by \cite{Hof}  \eqref{4.2} has at least a saturated equilibrium $x^*\ge 0$. The linearization of \eqref{4.2} at 0 is given by
$$x_i'(t)=\be_ix_i(t)+\sum_{j\ne i} d_{ij}x_j(t) , \quad i=1,2,\ldots,n. 
$$
With  $s(M(0))>0$,  this linear system is unstable (cf.~Theorem \ref{thm21}), and therefore the equilibrium 0  is not saturated, hence $x^*\ne 0$. If in addition $M(0)$ is irreducible, 0 is the only equilibrium of  \eqref{4.2} on the boundary of $\R^{n}_+$: otherwise (after a permutation of variables) there is an equilibrium of the form  $x^*=(0,\dots,0,x_{k+1}^*,\dots ,x_n^*)$ for some $k\in\{1,\dots,n-1\}$, then $\sum_{j=k+1}^n d_{ij}x_j^*=0$ for $1\le i\le k$, hence $d_{ij}=0$ for $1\le i\le k, k+1\le j\le n$, and $[d_{ij}]$ is not irreducible.
Therefore, we conclude that there  is an equilibrium of \eqref{4.2}  in the interior of $\R^n_+$. \end{pf}

\begin{rmk}\label{rmk41} Through the remainder of this section, for the matrix $M(0)$ in \eqref{M(0)}  we shall assume that  there is some positive vector $v$ such that $M(0)v>0$.
In the case of $M(0)$ an irreducible matrix, this condition can be simply replaced by the assumption  $s(M(0))>0$ (cf. Lemma \ref{lem3}).
\end{rmk}

By comparison with the cooperative system \eqref{2.4}, clearly Theorem \ref{thm32} provides an immediate criterion for  dissipativeness.

\begin{thm}\label{thm41} Suppose that $M(0)v>0$  for some positive vector $v$,
and let $X^*=(X_1^*,\dots,X_n^*)$ be a positive equilibrium for \eqref{2.4}, whose existence is given by Theorem \ref{thm31}.
Assume that $N_0$ is a non-singular M-matrix, and that $M(0)X^*>0$. Then,
system \eqref{eq} is dissipative; to be more precise, all positive solutions $x(t)$ of \eqref{eq} satisfy
\begin{equation}\label{4.1}
\limsup_{t\to\infty} x_i(t)\le X_i^*,\q i=1,\dots,n.
\end{equation}
If in addition  $M(0)$ is irreducible, system \eqref{2.1} has a positive equilibrium $x^*$.
\end{thm}


Next, we study the persistence of \eqref{eq}.
The  notion of persistence in Section 3 means that the population persists on each patch. We start with the  discussion  of  persistence  of the total population, therefore we refer to the more general concept of $\rho$-persistence as   in the monograph of Smith and Thieme \cite{ST}. Namely, with  $\rho(\var)=\sum_{i=1}^n \var_i(0)$, $\rho$-persistence means persistence of the total population.

%
%
%
%
%
%
%

\begin{thm}\label{thm42} 
Assume that \eqref{eq} is dissipative. If  $M(0)v>0$  for some positive vector $v$, then the total population is (weakly) persistent, i.e., 
$$\limsup_{t\to\infty}\sum_{i=1}^n x_i(t)>0$$
for all positive solutions $x(t)$ of \eqref{eq}.
Furthermore, if  $\be_i>0$ for $i=1,\dots,n$, then the total population is (strongly) uniformly persistent; i.e., there exists $\th>0$ such that
$$\liminf_{t\to\infty}\sum_{i=1}^n x_i(t)>\th $$
for all positive solutions $x(t)$ of \eqref{eq}.
\end{thm}

\begin{pf} Let $x(t)$ be a solution of  \eqref{eq}.  Since the system is dissipative,
$$\bar x_i:=\limsup_{t\to\infty} x_i(t)<\infty,\ \ \ i=1,\dots, n.$$
Choose $i\in\{ 1,\dots,n\}$ such that $\bar x_i=\max_{1\le j\le n} \bar x_j$. 
We  first claim that $\bar x_i>0$.

If $\bar x_i=0$, then $x_j(t)\to 0$ as $t\to\infty$ for all components $j$.
Take a positive vector $v$ such that $M(0)v>0$, and  
  choose $\vare>0$ small enough so that $(M(0)-\vare I)v>0$.
From Lemma \ref{lem1}, if $t$ is sufficiently large we have 
$$x_i'(t)\ge x_i(t)[\be_i -\vare -\mu_ix_i(t)]+\sum_{j\ne i} d_{ij}x_j(t-\tau_{ij}),\q i=1,\dots,n.$$
From Theorem \ref{thm32}, the cooperative system 
$$u_i'(t)= u_i(t)[\be_i-\vare -\mu_iu_i(t)]+\sum_{j\ne i} d_{ij}u_j(t-\tau_{ij}),\q i=1,\dots,n,$$ 
has a globally asymptotically stable equilibrium $u^*>0$. By comparison results \cite{Smith},
 we now obtain $\liminf_{t\to\infty}x_i(t)\ge u_i^*>0$, which is not possible. Therefore,  $\bar x_i>0$.
 \smal
 
Now, suppose that $\be_i>0$ for all $i$. By the fluctuation lemma there exists a sequence $(t_k)$ with $t_k\to \infty$, $x_i(t_k)\to \bar x_i$  and $ x_i'(t_k)\to 0$. Again from Lemma \ref{lem1}, for any $\vare >0$, if $k$ is sufficiently large we obtain
$$
x_i'(t_k)\ge x_i(t_k)\big [\be_i-\mu_ix_i(t_k)-(\bar x_i+\vare) \sum a_{ij}^+ \Big],\\
$$
By letting $k\to\infty$ and $\vare\to 0^+$, we obtain
$\bar x_i\ge \frac {\be_i}{\mu_i+\sum a_{ij}^+}>0.$
These arguments also show that
$$\limsup_{t\to\infty}\sum_{j=1}^n x_j(t)\ge \min_{1\le i\le n} \frac {\be_i}{\mu_i+\sum a_{ij}^+}=:\th_1>0.$$
Note that  the lower bound $\th_1$  does not depend on the particular solution $x(t)$. 
This means that the total population $\sum_{j=1}^n x_j(t)$ is {\it uniformly weakly persistent} (see \cite{ST} for a definition). On the other hand, since \eqref{eq} is dissipative, it has a compact global attractor \cite{Hale:1988}, and the hypotheses of Theorem 4.5 of \cite{ST}   are satisfied. This allows  to conclude the strong uniform persistence of the total population.

\end{pf}

\begin{cor}\label{cor41} 
Assume that  $\be_i>0$ for all $i$, and that  $N_0$ is a non-singular M-matrix. Then   \eqref{eq} is dissipative, the total population uniformly persists, and there exists a non-trivial  equilibrium $x^*\ge 0$. 
\end{cor}

Conditions for the  persistence of the population on each patch are given below. To simplify the notation,  denote
$$\widehat{N}=diag\, (\mu_1,\dots,\mu_n)-\Big [|a_{ij}|\Big].$$

\begin{thm}\label{thm43} 
Assume that  $M(0)v>0$  for some positive vector $v$ and that  $N_0$ is a non-singular M-matrix. If in addition $\widehat{N} X^*>0$,
where $X^*$ is the positive equilibrium of \eqref{2.4}, then \eqref{eq} is persistent and there is a positive equilibrium.
\end{thm}

\begin{pf} The existence of $X^*$, the unique positive equilibrium of \eqref{2.4}, is guaranteed by Theorem \ref{thm31}. Condition $\widehat{N} X^*>0$ translates as
\begin{equation}\label{4.3}
\mu_i X_i^*-\sum_j |a_{ij}|X_j^*>0,\ \ i=1,\dots,n,
\end{equation}
and in particular implies that $M(0)X^*=X^*\otimes N_0X^*>0$, i.e., 
$$
\mu_iX_i^*-\sum_j a_{ij}^- X^*_j>0,\q i=1,\dots,n.
$$
Theorem \ref{thm41} provides the upper bounds
$\bar x_i:=\limsup_{t\to\infty} x_i(t)\le X_i^*,$
for all $i$ and all solutions $x(t)$ of  \eqref{eq}. 
Now, define the matrix 
$$\widetilde{M}(0)=diag\, (\gamma_1,\dots,\gamma_n)+[d_{ij}],$$
where $ \gamma_i=\be_i-\sum_j a_{ij}^+X_j^*\, , 1\le i\le n.$
From \eqref{4.3}, 
\begin{equation}\label{4.4}
\begin{array}{ll}
(\widetilde{M}(0)X^*)_i&=X_i^*\Big [\be_i-\sum_j a_{ij}^+X_j^*\Big ]+\sum_{j\ne i}d_{ij}X_j^*\\
&=X_i^*\Big [ \mu_i X_i^*-\sum_j |a_{ij}|X_j^*\Big ]>0,\q 1\le i\le n.
\end{array}
\end{equation}

From Lemma \ref{lem1},  for each $\vare >0$  there exists $t_0>0$ such that 
$\int_0^{\infty}K_{ij}(s)x_j(t-s)\, ds\le (1+\vare)X_j^*$ for any $i,j=1,\dots,n$ and $t\ge t_0$. Thus, for $t\ge t_0$,
 $$
\begin{array}{ll}
\displaystyle{x_i'(t)\ge x_i(t)\bigg(\be_i-(1+\vare)\sum_{j} a_{ij}^+X_j^*-\mu_i x_i(t)+\sum_{j} a_{ij}^- \int_0^{\infty}K_{ij}(s)x_j(t-s)ds \bigg)} \\
\displaystyle{\hskip 1cm +\sum_{j\ne i} d_{ij}x_j(t-\tau_{ij}), \q i=1,2,\ldots,n}.
\end{array}
$$
In virtue of \eqref{4.4}, we can choose $\vare>0$ small enough so that 
$\widetilde{M}_\vare (0)X^*>0$, where
$$\widetilde{M}^\vare (0)= diag\, (\gamma_1^\vare,\dots,\gamma_n^\vare)+[d_{ij}],\q {\rm for }\q \gamma_i^\vare=\be_i-(1+\vare)\sum_j a_{ij}^+X_j^*\, .$$
From Theorem \ref{thm31},  observe that  the cooperative system
\begin{equation}\label{4.5}
\begin{array}{ll}
\displaystyle{u_i'(t)=u_i(t)\bigg (\gamma_i^\vare-\mu_i u_i(t)+\sum_{j} a_{ij}^- \int_0^{\infty}K_{ij}(s)u_j(t-s)ds \bigg)+\sum_{j\ne i} d_{ij}u_j(t-\tau_{ij}),\  i=1,2,\ldots,n},
\end{array}
\end{equation}
is persistent.
Comparing the solutions of \eqref{eq} with the solutions of  \eqref{4.5}, we deduce that \eqref{eq} is persistent as well. Now, from the persistence and Theorem \ref{thm41}, there is a positive equilibrium.
\end{pf}

We finally present a criterion for the global asymptotic stability of a positive equilibrium for \eqref{eq}.

\begin{thm}\label{thm44} 
Assume that \eqref{eq} is dissipative, persistent and has  an equilibrium $x^*>0$.
If in addition $\widehat{N}x^*>0$, then  $x^*$ is globally attractive. 
\end{thm}

\begin{pf}  By the change of variables $y(t)=x(t)-x^*$,  \eqref{eq} becomes
\begin{equation}\label{4.8}
\begin{array}{ll}
\displaystyle{y_i'(t)=-(y_i(t)+x_i^*)\biggl(\mu_i y_i(t)+\sum_j a_{ij} \int_0^{\infty}K_{ij}(s)y_j(t-s)ds\biggr)} 
\\
\displaystyle{\hskip 2cm -y_i(t)\frac{1}{x_i^*}\sum_{j\ne i} d_{ij}x_j^*+\sum_{j\ne i} d_{ij}y_j(t-\tau_{ij}), \quad i=1,2,\ldots,n}.
\end{array}
\end{equation}

Define $u_j=\limsup_{t\to \infty} y_j(t), -v_j=\liminf_{t\to \infty} y_j(t)$, and
$$U=\max_j\frac{u_j}{x_j^*},\q V=\max_j\frac{v_j}{x_j^*},\q L=\max (U,V).$$
Clearly $L\ge 0$. Moreover, from the persistence $v_j<x_j^*$ for all components $j$, and therefore $V<1$.

It suffices  to show that $L=0$. 
We argue by contradiction, so  assume $L>0$. 

Consider first the case of $L=U$, and choose $i$ such that $U=\frac{u_i}{x_i^*}$. Take a sequence $t_k\to \infty$ such that $y_i(t_k)\to u_i$ and $y_i'(t_k)\to 0$. Applying Lemma \ref{lem1} to \eqref{4.8}, for any $\vare >0$, if $k$ is sufficiently large we obtain
$$
\begin{array}{ll}
\displaystyle{y_i'(t_k)\le -(y_i(t_k)+x_i^*)\biggl(\mu_i y_i(t_k)-(1+\vare)L\sum_j |a_{ij}|x_j^*\biggr) -y_i(t_k)\frac{1}{x_i^*}\sum_{j\ne i} d_{ij}x_j^*} \\
\displaystyle{\hskip 3.5cm+(1+\vare)L\sum_{j\ne i,} d_{ij}x_j^*}.
\end{array}
$$
By letting $k\to\infty$ and $\vare\to 0^+$, we get\begin{equation}\label{4.10}
\begin{array}{l}
\displaystyle{0\le-(L+1)x_i^*\biggl(\mu_i x_i^*-\sum_{j=1}^n |a_{ij}| x_j^*\biggr) L}<0,
\end{array}
\end{equation}
a contradiction.
Now, consider the case $L=V=\frac{v_i}{x_i^*}$ for some $i$. Then, there is a sequence $t_k\to \infty$ with  $y_i(t_k)\to -v_i=-Lx_i^*>-x_i^*$ and $y_i'(t_k)\to 0$. 
We proceed as  in the  above case and   instead of \eqref{4.10}  obtain
\begin{equation}
 \begin{array}{l}
\displaystyle{0\ge(-L+1)x_i^*\biggl(\mu_i x_i^*-\sum_{j=1}^n |a_{ij}| x_j^*\biggr) L}>0,
\end{array}
\end{equation}  which is again a contradiction. The proof is complete.
\end{pf}

By Theorems \ref{thm41},  \ref{thm43} and   \ref{thm44}, we immediately get:

\begin{cor}\label{cor42} 
Assume that  $M(0)v>0$  for some positive vector $v$,   $N_0$ is a non-singular M-matrix, and $\widehat{N}X^*>0$, where $X^*$ is the positive equilibrium of \eqref{2.4}. Then there exists an equilibrium $x^*>0$ of \eqref{eq}.
If in addition $\widehat{N}x^*>0$, then  $x^*$ is globally attractive. 
\end{cor}


If $M(0)$ is a positive matrix, the assumption $\widehat{N}X^*>0$ can be dropped in the above criterion, since  one can use the persistence of the total population, rather than the persistence on each patch.

\begin{thm}\label{thm45} 
Assume that  $N_0$ is a non-singular M-matrix, and $\be_i,d_{ij}>0$ for all $i,j=1,\dots,n$.
If  the  equilibrium $x^*>0$ of \eqref{eq} (whose existence is given in Theorem \ref{thm41}) satisfies $\widehat{N}x^*>0$, then  $x^*$ is globally attractive.  
\end{thm}

\begin{pf}  Under the  assumption $\be_i,d_{ij}>0$ for all $i,j$, by Theorems \ref{thm41} and \ref{thm42}, system \eqref{eq} is dissipative, the total population is uniformly persistent, and there is an equlibrium $x^*>0$.

We now use the same notation and proceed as  in the proof of Theorem \ref{thm44}, noting however that $V\le 1$, but the situation $V=1$ is possible. In fact, $v_j\le x_j^*$ for all $j$, and $v_j<x_j^*$ for at least one component $j$, because of the persistence of the total population. By repeating that proof, we only have to further  assure that the  case of $L=V=1$ is not possible.

Let $L=V=\frac{v_i}{x_i^*}=1$ for some $i$. 
Consider a sequence $t_k\to \infty$ with  $y_i(t_k)\to -v_i=-x_i^*$ and $y_i'(t_k)\to 0$.  Applying Lemma \ref{lem1} to \eqref{4.8}, for any $\vare >0$, if $k$ is sufficiently large we obtain
$$
\begin{array}{ll}
\displaystyle{y_i'(t_k)\ge -(y_i(t_k)+x_i^*)\biggl(\mu_i y_i(t_k)+(1+\vare)\sum_j |a_{ij}|x_j^*\biggr)} \\
\displaystyle{\hskip 2cm -y_i(t_k)\frac{1}{x_i^*}\sum_{j\ne i} d_{ij}x_j^*+(1+\vare)\sum_{j\ne i} d_{ij}v_j}.
\end{array}
$$
By letting $k\to\infty$ and $\vare\to 0^+$, we obtain  a contradiction, since
$
0\ge -\sum_{j\ne i} d_{ij}x_j^*+\sum_{j\ne i} d_{ij}v_j>0.
$
\end{pf}

\begin{exmp} Consider the following system of the form \eqref{eq} with $n=2$:
\begin{equation}\label{4.14}
\begin{array}{ll}
\dps{x_1'(t)=x_1(t)\biggl(\be_1-\mu_1x_1(t)- a_{11} \int_0^{\infty}K_{11}(s)x_1(t-s)\, ds} \\
\dps{\hskip 2cm - a_{12} \int_0^{\infty}K_{12}(s)x_2(t-s)\, ds\biggr)+d_{1}x_2(t-\tau_1)}\\
\dps{x_2'(t)=x_2(t)\biggl(\be_2-\mu_2x_2(t)- a_{21} \int_0^{\infty}K_{21}(s)x_1(t-s)\, ds} \\
\dps{\hskip 2cm - a_{22} \int_0^{\infty}K_{22}(s)x_2(t-s)\, ds\biggr)+d_{2}x_1(t-\tau_2)}.
\end{array}
\end{equation}
with  delays $\tau_1,\tau_2\ge 0$ and coefficients $d_i>0, a_{ij}\ge 0$ and $\be_i\in\R, i,j=1,2$. For this system, and with the previous notation,  
$$
M(0)=\left[
\begin{array}{cc}
\be_1&d_1 \\
d_2& \be_2 \\
\end{array}
\right],\ N_0=\left[
\begin{array}{cc}
\mu_1&0 \\
0& \mu_2 \\
\end{array}
\right].$$
Note that $M(0)$ is irreducible. 
We have $s(M(0))\le 0$   if and only if $\be_1\le 0,\be_2\le 0$ and $\be_1\be_2\ge d_1d_2$, in which case the trivial equilibrium is a global attractor of all positive solutions (cf. Theorem \ref{thm22}); otherwise, Theorem \ref{thm41} assures that there exists a positive equilibrium. 

As an illustration, now take \eqref{4.14} subject to the constraints
\begin{equation}\label{4.15}
\begin{array}{ll}
 \mu_1>a_{11}+a_{12}, \ \mu_2>a_{12}+a_{22}, \\ 
 \\
\dps\frac{\be_1+d_1}{\mu_1+a_{11}+a_{12}}=\frac{\be_2+d_2}{\mu_2+a_{12}+a_{22}}=:c>0. \\
\end{array}
\end{equation}
Under these conditions, one easily verifies that $s(M(0))>0$, and that $x^*=(c,c)$  is an equilibria of  \eqref{4.14}.
 The matrix $\widehat{N}$ reads as
$\widehat{N}=\left[
\begin{array}{cc}
\mu_1-a_{11}&-a_{12} \\
-a_{21}& \mu_2-a_{22} \\
\end{array}
\right],$
so $\widehat{N}x^*>0$. Hence, if $\be_1>0,\be_2>0$,  Theorem \ref{thm45} implies that  $x^*=(c,c)$ is a global attractor of all positive solutions of \eqref{4.14}.

 For the situation  $\be_1\le 0$ or $\be_2\le 0$,  together with \eqref{4.15} if we now assume
$$\dps \frac{\be_1+d_1}{ \mu_1}=\frac{\be_2+d_2}{\mu_2}=:\gamma,$$
then $X^*=(\gamma,\gamma)$ is a globally attractive  equilibrium
for the  cooperative system associated with \eqref{4.14}:
\begin{equation}\label{4.16}
\begin{array}{ll}
\dps{x_1'(t)=x_1(t) (\be_1-\mu_1 x_1(t))+d_{1}x_2(t-\tau_1)}\\
\dps{x_2'(t)=x_2(t)(\be_2-\mu_2 x_2(t))+d_{2}x_1(t-\tau_2)}.
\end{array}
\end{equation}
 Since $\widehat{N}X^*>0$, from Corollary \ref{cor42} then  $x^*=(c,c)$  globally attracts the  positive solutions of \eqref{4.14}. 
\end{exmp}


	\begin{rmk} The requirements  $\widehat{N}X^*>0$ and $\widehat{N}x^*>0$ in Theorems \ref{thm43}, \ref{thm44} and  \ref{thm45} are expressed in terms of the positive equilibria $X^*$ of \eqref{2.4} and   $x^*$ of \eqref{eq}. It would be therefore relevant to improve the above criteria, in the sense of achieving  sufficient conditions for the uniform persistence  of \eqref{eq} and the global attractivity of $x^*$ involving only the coefficients of the system. 
The theorem below is a first attempt to establish such type of criteria.	
\end{rmk}

\begin{thm}\label{thm46} Suppose that, for all $i,j=1,\dots,n$,
\begin{equation}\label{4.21}
 \begin{array}{ll}
& {\mu_i >\sum_j a_{ij}^-,}\\
&  {\be_i\ge Ma_{ii}^+,\q d_{ij}\ge M a_{ij}^+,\q i\ne j}\\
 \end{array}
\end{equation}
with $\be_i+\sum_{j\ne i}d_{ij}>M \sum_ja_{ij}^+$, where
\begin{equation}\label{4.7''}
M=\max_{1\le i\le n} \frac{\be_i+\sum_{j\ne i}d_{ij}}{\mu_i-\sum_j a_{ij}^-}.
\end{equation}
Then  \eqref{eq} is dissipative and persistent. If, for all $i,j=1,\dots,n$,
\begin{equation}\label{4.7}
 \begin{array}{ll}
& {\mu_i >\sum_j a_{ij}^-,}\\
&  {\be_i\ge 2Ma_{ii}^+,\q d_{ij}\ge 2M a_{ij}^+,\q i\ne j}\\
 \end{array}
\end{equation}
with $\be_i+\sum_{j\ne i}d_{ij}>2M \sum_ja_{ij}^+$,
then  \eqref{eq} has an equilibrium $x^*>0$ which  is  globally attractive.  
\end{thm}

\begin{pf} If \eqref{4.21} holds, we have $M(0)q>0$ and $N_0q>0$, for $q=(1,\dots,1)$. From Theorem \ref{thm32}, we derive that there exists a positive equilibrium $X^*=(X_1^*,\dots,X_n^*)$ of \eqref{2.4}, with $M(0)X^*>0$. For $X_i^*=\max_{1\le j\le n} X_j^*$, one easily checks that $X_i^*\le \frac{\be_i+\sum_{j\ne i}d_{ij}}{\mu_i-\sum_j a_{ij}^-}$, and hence  the estimates $X_j^*\le M$ 
for $M$ defined in \eqref{4.7''}. 
 To conclude the persistence of \eqref{eq}, from Theorem \ref{thm43} it is sufficient to show that $\widehat{N}X^*>0$. From the identities 
$$
\be_iX_i^*+\sum_{j\ne i}d_{ij}X_j^*=X_i^*(\mu_iX_i^*-\sum_j a_{ij}^-X_j^*),\q i=1,\dots,n,$$
 we deduce that 
$$
\begin{array}{ll}
X_i^*(\widehat{N} X_i^*)_i&=
X_i^*\Big (\mu_i X_i^*-\sum_j |a_{ij}|X_j^*\Big )\\
&=X_i^*\Big (\mu_i X_i^*-\sum_j a_{ij}^-X_j^*-\sum_j a_{ij}^+ X_j^*\Big)\\
&=(\be_i-a_{ii}^+X_i^*)X_i^*+\sum_{j\ne i}(d_{ij}-a_{ij}^+X_i^* ) X_j^*\\
&\ge (\be_i-Ma_{ii}^+)X_i^*+\big (\min\limits_{1\le j\le n} X_j^* \big)\sum_{j\ne i}\big (d_{ij}-Ma_{ij}^+\big) > 0.
\end{array}
$$

Next, suppose that the stronger conditions \eqref{4.7} hold. The  components of the positive equilibrium $x^*=(x_1^*,\dots ,x_n^*)$ of
 \eqref{eq} also satisfy the estimates
 $x_i^*\le M,  i=1,\dots,n,$
 for $M$ as in \eqref{4.7''}. 
Proceeding as above, $\widehat{N}x^*>0$, hence the conclusion follows from  Theorem \ref{thm45}. Details are omitted.
\end{pf}

\begin{exmp} Consider again the system \eqref{4.14}, with all coefficients being positive. If $\mu_1\ge \be_1+d_1, \mu_2\ge \be_2+d_2, \be_1>2a_{11}, \be_2>2a_{22}, d_1\ge 2a_{12}, d_2\ge 2a_{21},$
then $M\le 1$ for $M$ as in \eqref{4.7''} and the constraints \eqref{4.7} are fulfilled, hence there is a positive equilibrium which is a global attractor.
\end{exmp}

\section*{Acknowledgement} The research was supported by Funda\c c\~ao para a Ci\^encia e a Tecnologia (Portugal),  PEst-OE/\-MAT/\-UI0209/2011.

\end{document}